\documentclass[leqno,10pt]{article}
 
\usepackage{amsmath,amsthm}
\usepackage{amssymb}
\usepackage{euscript}
\usepackage[frame,ps,matrix,arrow,curve,rotate]{xy} 

\numberwithin{equation}{subsection}  % [if desired]

\newcommand{\sqsp}{\renewcommand{\baselinestretch}{1.1}\tiny\normalsize}

\newtheorem{thm}[equation]{Theorem}
\newtheorem{prop}[equation]{Proposition}
\newtheorem{cor}[equation]{Corollary}
\newtheorem{lemma}[equation]{Lemma}

\newtheorem{thm-intro}{Theorem}

\theoremstyle{definition}

\newcommand{\bC}{\mathbf{C}}

\newcommand{\bQ}{\mathbf{Q}}

\newcommand{\bW}{\mathbf{W}}
\newcommand{\bZ}{\mathbf{Z}}

\DeclareMathOperator{\Id}{Id} 
\DeclareMathOperator{\Ring}{\mathbf{Ring}}      % the category of rings

\newcommand{\llbrack}{\lbrack \lbrack}
\newcommand{\rrbrack}{\rbrack \rbrack}

\begin{document}
\title{Moduli space of filtered lambda-ring structures over a filtered ring}
\author{Donald Yau}

\maketitle

\begin{abstract}
Motivated by recent works on the genus of classifying spaces of compact Lie groups, here we study the set of filtered $\lambda$-ring structures over a filtered ring from a purely algebraic point of view.  From a global perspective, we first show that this set has a canonical topology compatible with the filtration on the given filtered ring.  For power series rings $R \llbrack x \rrbrack$, where $R$ is between $\bZ$ and $\bQ$, with the $x$-adic filtration, we mimic the construction of the Lazard ring in formal group theory and show that the set of filtered $\lambda$-ring structures over $R \llbrack x \rrbrack$ is canonically isomorphic to the set of ring maps from some ``universal'' ring $U$ to $R$.  From a local perspective, we demonstrate the existence of uncountably many mutually non-isomorphic filtered $\lambda$-ring structures over some filtered rings, including rings of dual numbers over binomial domains, (truncated) polynomial and powers series rings over torsionfree $\bQ$-algebras.
\end{abstract}

\noindent 
E-mail: \texttt{dyau@math.uiuc.edu} \\
\noindent
Address: Department of Mathematics, University of Illinois at Urbana-Champaign, 1409 W.\ Green Street, Urbana, IL 61801, USA \\

\noindent
\underline{\hspace{5in}}

%\tableofcontents

\sqsp
                                   %%  Introduction  %%

\section{Introduction}
\label{sec:introduction}
Let $R$ be a ring, by which we always mean a commutative ring with unit.  A $\lambda$-ring (originally called \emph{special} $\lambda$-ring) structure on $R$ in the sense of Grothendieck consists of functions $\lambda^i \colon R \to R$ ($i \geq 0$) satisfying properties analogous to those satisfied by exterior power operations.  Examples of $\lambda$-rings include the complex representation ring $R(G)$ of a group $G$, the complex $K$-theory $K(X)$ of a space $X$, and the equivariant $K$-theory $K_G(X)$ of a $G$-space $X$.  A filtered ring is a ring $R$ which comes equipped with a multiplicative decreasing filtration $\lbrace I^n \rbrace$ of ideals.  A filtered $\lambda$-ring is a filtered ring $R$ which is also a $\lambda$-ring for which the filtration ideals are all closed under the operations $\lambda^i$ for $i > 0$.

The main aim of this paper is to study the set of filtered $\lambda$-ring structures over a filtered ring.  Before discussing the results of this paper, let us first provide some motivations for this work, which comes from the author's recent study of the genus of classifying spaces of compact connected Lie groups.  

Recall that the \emph{genus} of a nilpotent finite type space $X$ consists of the homotopy types of all nilpotent finite type spaces $Y$ such that the $p$-completions of $X$ and $Y$ are homotopy equivalent for all primes $p$ and also their rationalizations are homotopy equivalent.  The genus of a classifying space is a very rich set.  Indeed, M$\o$ller showed in \cite{moller} that the genus of the classifying space $BG$ of any non-abelian connected compact Lie group $G$ is an uncountable set.  The case when $G$ is $SU(2)$ was first studied by Rector \cite{rector}.  To understand these uncountably many homotopically distinct spaces, one would like to find an algebraic invariant which has the properties that (1) it can tell apart different spaces in the genus of $BG$, and that (2) it is computable to some extent.  
An algebraic invariant with the first property was found by Notbohm in \cite{notbohm}, in which he showed that $K$-theory $\lambda$-rings classify the genus of $BG$ provided $G$ is a simply-connected compact Lie group.  The question, therefore, became how to compute these $K$-theory $\lambda$-rings and how to use them to answer topological questions.

The author showed in \cite{yau1} that $K$-theory \textit{filtered} ring is an invariant of the genus of a space $X$, provided $X$ is simply-connected and has even and torsionfree integral homology and a finitely generated power series $K$-theory filtered ring.  Here the filtration in $K$-theory comes from a skeletal filtration of the space.  More precisely, if $X$ is such a space and $Y$ is in its genus, then there exists a filtered ring isomorphism from $K(X)$ to $K(Y)$.  This applies to, for example, $BSp(n)$ for any $n \geq 1$, $BSU(n)$ for any $n \geq 2$, or any finite product of copies of such spaces and of infinite complex projective space.

An immediate, purely algebraic consequence of these results, when applied to the spaces $BSp(n)$, is the following:  For any positive integer $n$, the filtered power series ring $\bZ \llbrack x_1, \ldots, x_n \rrbrack$ (the $K$-theory ring of $BSp(n)$) on $n$ indeterminates, where $x_i$ has filtration precisely $4i$, admits uncountably many mutually non-isomorphic filtered $\lambda$-ring structures.

These results just described also have significant topological implications.  In fact, by combining it with a theorem of Atiyah \cite{atiyah}, the author \cite{yau2} was able to completely classify spaces in the genus of $BSU(2)$ which admit non-nullhomotopic maps from infinite complex projective space (``the maximal torus") and to compute the maps explicitly in these cases.  This extends the classical result of McGibbon \cite{mcgibbon} and Rector \cite{rector} on maximal torus admissibility of spaces in the genus of $BSU(2)$.

Of course, different filtered $\lambda$-ring structures over an underlying filtered ring are already interesting from a purely algebraic point of view.  But with the above topological results in mind, it becomes clear that they should also be studied because of their potential to provide information about spaces in the same genus.  
With these motivations in mind, we now describe the main results of this paper.

In the first two sections below we will take on a \emph{global} perspective, in the sense that we consider all the filtered $\lambda$-ring structures over a given filtered ring $R$ simultaneously.  We first address the following moduli problem:

\begin{quote}
\textit{Is there a canonical way to put a topology on the set of filtered $\lambda$-ring structures over a given filtered ring $R$, which is in some sense compatible with the filtration on R?}
\end{quote}
The answer to this question involves the universal ring of Witt vectors $\bW(R)$, sometimes called the big Witt vectors on $R$.  The functor $\bW$ on the category of rings is actually (part of) a comonad.  The following theorem, which gives an affirmative answer to the above question, is an amalgamation of the results in \S \ref{sec:moduli}.

\bigskip
\begin{thm-intro}
\label{thm-intro:moduli}
The following statements hold.
\begin{itemize}
\item $\bW$ restricts to a functor on the category of filtered rings.
\item $\bW$ is $($part of$)$ a comonad on the category of filtered rings.
\end{itemize}
Now regard $\bW$ as a comonad on the category of filtered rings and let $R$ be a filtered ring.  Then:
\begin{itemize}
\item $\bW(R)$ has a canonical filtered $\lambda$-ring structure.
\item $R$ is complete $($resp.\ Hausdorff$)$ if and only if $\bW(R)$ is so.
\item There is a canonical bijection between the set of filtered $\lambda$-ring structures on $R$ and the set of $\bW$-coalgebra structures on $R$.  In particular, the former has a canonical topology as the subspace of the space of continuous maps from $R$ to $\bW(R)$ consisting of the filtered ring homomorphisms which make $R$ into a $\bW$-coalgebra.
\end{itemize}
\end{thm-intro}

We will refer to this as the \emph{moduli space of filtered} $\lambda$-\emph{ring structures on} $R$.

The proof of this theorem consists of extending known results in the non-filtered setting to the filtered one.  In fact, it is a well-known fact that for any ring $R$, a $\lambda$-ring structure on $R$ corresponds precisely to a $\bW$-coalgebra structure on $R$, where $\bW$ is regarded as a comonad on the category of rings.  So what we will do is to show that this correspondence restricts to the (non-full) subcategory of filtered rings and filtered ring homomorphisms.  As in the non-filtered case, it is more convenient to work with a functor $\Lambda$ which is naturally isomorphic to $\bW$ via the so-called exponential isomorphism $E$.

Now one way to obtain filtered rings is to consider power series rings $R \llbrack x \rrbrack$ with the $x$-adic filtration, where $R$ is any ring and $x$ is an indeterminate with a fixed positive filtration.  For example, the $K$-theory of $BSU(2)$ is the filtered power series ring $\bZ \llbrack x \rrbrack$ with $x$ in filtration precisely $4$.  It is, therefore, natural to ask the following representability question:

\begin{quote}
\textit{Is the functor which sends a ring $R$ to the set of filtered $\lambda$-ring structures on $R \llbrack x \rrbrack$ $($co$)$representable?}
\end{quote}
The reader will notice that this question is similar to one in formal group theory: Is the functor which sends a ring $R$ to the set of formal group laws over it (co)representable?  The affirmative answer to this formal group question is a well-known result of Lazard, which says that this functor is co-represented by what is now called the Lazard ring $L$.  The isomorphism from $\Ring(L,R)$ (the set of ring maps from $L$ to $R$) to the set of formal group laws over $R$ is given as follows:  It sends a ring map $f \colon L \to R$ to the formal group law over $R$ which is the push-forward by $f$ of the universal formal group law over $L$.  We will approach our representability problem in a similar way.

In order to exploit the analogy with formal group theory, we would like to reduce a filtered $\lambda$-ring structure over $R \llbrack x \rrbrack$ to a sequence of power series $\psi^n(x)$ over $R$ (the Adams operations) satisfying the properties: 
\begin{enumerate}
\item $\psi^1(x) = x$,
\item $\psi^m (\psi^n(x)) = \psi^n(\psi^m(x)) = \psi^{mn}(x)$,
\item $\psi^p(x) \equiv x^p$ (mod $pR$) for any prime $p$.
\end{enumerate}
As we will see in \S \ref{sec:power series} below, this is possible provided that we impose the heavy restriction, $\bZ \subseteq R \subseteq \bQ$.

To state our next result, recall that a \emph{binomial domain} is a $\bZ$-torsionfree ring $R$ for which the binomial symbol
\[
\binom{r}{n} = \frac{r(r-1) \cdots (r-n+1)}{n!} \in R \otimes \bQ
\]
lies in $R$ for any element $r$ in $R$ and positive integer $n$.  For example, $\bZ$, $\bZ_{(p)}$ ($p$-local integers), $\bZ_p$ ($p$-adic integers), and torsionfree $\bQ$-algebras (and hence fields of characteristic $0$) are binomial domains.

We then have the following theorem, which is a simplified version of the main result (Theorem \ref{thm:co-representing ring}) in \S \ref{sec:power series}.

\bigskip
\begin{thm-intro}
\label{thm-intro:Lazard}
There exists a binomial domain $U$ such that for every ring $R$ between $\bZ$ and $\bQ$, there is a canonical isomorphism
\[
\Psi_R \colon \Ring(U,R) ~\xrightarrow{\cong}~ \left\lbrace\text{filtered }\lambda \text{-ring structures over } R\llbrack x \rrbrack\right\rbrace.
\]
\end{thm-intro}

The method of proof of this theorem has some similarity with the construction of the Lazard ring.  Roughly speaking, one adjoins to $\bZ$ variables corresponding to the coefficients of the power series $\psi^p(x)$ and imposes appropriate relations, though the actual construction is slightly more complicated than this.

It should also be remarked that this ring $U$ is non-trivial, since the power series filtered ring $\bZ \llbrack x \rrbrack$ admits uncountably many mutually non-isomorphic filtered $\lambda$-ring structures (see \cite{yau1}).

The above two theorems are our main results about filtered $\lambda$-ring structures over a filtered ring from a global perspective.

In \S \ref{sec:examples} we take on a \emph{local} perspective and consider the points in the moduli spaces of filtered $\lambda$-ring structures over certain filtered rings.  This consists of counting the number of points and, therefore, distinguishing between points in the moduli spaces.  We will consider the following filtered rings:
\begin{description}
\item[Ring of dual numbers] $B \lbrack \varepsilon \rbrack$ ($\varepsilon^2 = 0$), where $B$ is a non-zero binomial domain,
\item[Possibly truncated polynomial rings] $A \lbrack x \rbrack$ and $A \lbrack x \rbrack/x^N$, where $N \geq 2$ and $A$ is a non-zero, $\bZ$-torsionfree $\bQ$-algebra,
\item[Power series ring] $A \llbrack x \rrbrack$, where $A$ is a non-zero, $\bZ$-torsionfree $\bQ$-algebra.
\end{description}
Each of these rings is equipped with the $x$-adic (or $\varepsilon$-adic) filtration.

Then we have the following theorem, which is extracted from \S \ref{sec:examples} (see Theorems \ref{thm:dual numbers} and \ref{thm:polynomial ring}).

\bigskip
\begin{thm-intro}
\label{thm-intro:examples}
Each of these three filtered rings admits uncountably many mutually non-isomorphic filtered $\lambda$-ring structures.
\end{thm-intro}

A few remarks are in order.  First, the conclusion in the theorem is stronger than just saying that there are uncountably many points in the moduli spaces.  Second, the proofs of this existence result for the various rings are similar.  In each case, we write down the Adams operations and use Wilkerson's Theorem \ref{thm:wilkerson} to conclude that they come from a $\lambda$-ring structure.  Third, for the ring of dual numbers over a non-zero binomial domain, we actually obtain a complete classification of all of the isomorphism classes of filtered $\lambda$-ring structures over it.

This finishes the presentation of our main results.  The rest of the paper is organized as follows.

In \S \ref{sec:moduli} we begin by recalling the definition of a $\lambda$-ring and some properties of the functor $\Lambda$ on the category of rings and prove Theorem \ref{thm-intro:moduli} in this setting.  Then we transfer the result to the functor $\bW$ by studying the exponential isomorphism.

In \S \ref{sec:power series} we first introduce the notion of a (filtered) $\psi$-ring, which one can think of as a (filtered) ring with Adams operations.  For a ring $R$ between $\bZ$ and $\bQ$, a filtered $\lambda$-ring over $R \llbrack x \rrbrack$ is precisely a filtered $\psi$-ring over it.  By mimicking the construction of the Lazard ring, we then construct the ring $U$ which comes equipped with a  ``universal" $\psi$-ring structure over $U \llbrack x \rrbrack$ and prove a more precise version of Theorem \ref{thm-intro:Lazard}.

Finally, Theorem \ref{thm-intro:examples} is proved in \S \ref{sec:examples}.  We will end this paper with the observation (Theorem \ref{thm:adams operations}) that, for any domain $R$ (that is, a commutative ring with unit without zero divisors), a filtered $\lambda$-ring structure on $R \llbrack x \rrbrack$ is determined by any single Adams operation $\psi^p$, $p$ prime.  This depends on a result (Proposition \ref{prop:lubin}) about commuting power series which might be of independent interest.

                                %%  Acknowledgement  %%
\section*{Acknowledgment}
I am indebted to Lars Hesselholt and Haynes Miller for sharing with me their knowledge and ideas about $\lambda$-rings in particular and mathematics in general.  Theorem \ref{thm-intro:moduli} and Theorem \ref{thm-intro:Lazard} are the results of questions of, respectively, Lars Hesselholt and Haynes Miller.  I also want to thank Kristine Baxter Bauer, Nitu Kitchloo, Peter May, and Randy McCarthy for discussions about this work.

                                  %%  Moduli space  %%

\section{Moduli space of filtered $\lambda$-ring structures}
\label{sec:moduli}

We remind the reader that all rings considered here are commutative with unit, and so by a ring we mean a commutative ring with unit.  

The purpose of this section is to prove Theorem \ref{thm-intro:moduli} in the Introduction.  We first recall the definition of a (filtered) $\lambda$-ring.

\subsection{Filtered $\lambda$-ring}
\label{subsec:filtered lambda-ring}

A \emph{filtered ring} is a ring $R$ equipped with a decreasing filtration 
\[
R ~=~ I^0 ~\supset~ I^1 ~\supset~ \cdots
\]
of ideals, called filtration ideals, such that $I^nI^m \subseteq I^{n+m}$ for any $n, m \geq 0$.  A filtered ring homomorphism $f \colon R \to S$ between two filtered rings is a ring map which preserves the filtrations: $f(I^n_R) \subset I^n_S$ for all $n$.  The category of filtered rings and filtered ring homomorphisms is regarded as a subcategory of the category of rings.  The filtration ideals are sometimes suppressed from the notation.

A $\lambda$-\emph{ring} is a ring $R$ equipped with functions $\lambda^i \colon R \to R$ $(i \geq 0)$, called $\lambda$-operations, such that for any elements $r$ and $s$ in $R$, the following conditions hold:
\begin{itemize}
\item $\lambda^0(r) = 1$,
\item $\lambda^1(r) = r$,
\item $\lambda^n(1) = 0$ for all $n > 1$,
\item $\lambda^n(r + s) = \sum_{i = 0}^n \lambda^i(r) \lambda^{n-i}(s)$,
\item $\lambda^n(rs) = P_n(\lambda^1(r),\ldots,\lambda^n(r);\lambda^1(s),\ldots,\lambda^n(s))$,
\item $\lambda^m(\lambda^n(r)) = P_{m,n}(\lambda^1(r),\ldots,\lambda^{mn}(r))$.
\end{itemize}
Here the $P_n$ and $P_{m,n}$ are certain universal polynomials with integer coefficients (see Atiyah-Tall \cite{at} or Knutson \cite{knutson} for detail).

A \emph{filtered} $\lambda$-\emph{ring} is a filtered ring $R$ which is also a $\lambda$-ring for which the filtration ideals are all closed under the operations $\lambda^i$ ($i > 0$).  A filtered $\lambda$-ring map $f \colon R \to S$ between two filtered $\lambda$-rings is a filtered ring map which commutes with the $\lambda$-operations: $f\lambda^i = \lambda^i f$ for all $i$.

%%%%%%%%%%%%%%%%%%%%%%%%%%%%%%%%%%%
%%%%%%%%%%%%%%%%%%%%%%%%%%%%%%%%%%%

\subsection{Comonads and coalgebras}
\label{subsec:comonads} 
Since several results later on in this section use the notions of a comonad and of an algebra over a comonad, we will briefly recall the relevant definitions here.  The reader can consult MacLane's book \cite{maclane} for more information on this topic.

A \emph{comonad} on a category $\bC$ is an endofunctor $F \colon \bC \to \bC$ together with natural transformations $\eta \colon F \to \Id$ and $\mu \colon F \to F^2$, called counit and comultiplication, such that the following counital and coassociativity conditions hold:
\[
F\eta \circ \mu ~=~ \Id_R ~=~ \eta F \circ \mu \qquad \text{and} \qquad
F\mu \circ \mu ~=~ \mu F \circ \mu.
\]
The natural transformations $\mu$ and $\eta$ are often omitted from the notation and we speak of $F$ as a comonad.

If $F$ is a comonad on a category $\bC$, then an $F$-{\it coalgebra structure} on an object $X$ of $\bC$ is a morphism $\xi \colon X \to FX$ in $\bC$ such that the following counital and coassociativity conditions hold:
\[
F\xi \circ \xi ~=~ \mu \circ \xi \qquad \text{and} \qquad \eta \circ \xi ~=~ \Id_X.
\]  
A map of $F$-coalgebras $g \colon (X, \xi_X) \to (Y, \xi_Y)$ consists of a morphism $g \colon X \to Y$ in $\bC$ such that $\xi_Y \circ g = Fg \circ \xi_X$.

%%%%%%%%%%%%%%%%%%%%%%%%%%%%%%%%%%%%%%%%%
%%%%%%%%%%%%%%%%%%%%%%%%%%%%%%%%%%%%%%%%%

\subsection{The functor $\Lambda$}
\label{subsec:Lambda}
Here we recall only those properties of $\Lambda$ which we will need later on in this section.  A general reference for the functor $\Lambda$ (in the category of rings) is Hazewinkel's book \cite[Ch.\ 3]{hazewinkel}.

If $A$ is any ring, then $\Lambda(A)$, sometimes called the universal $\lambda$-ring associated to $A$, is the ring which as a set consists of the power series $f(t) = 1 + \sum_{i \geq 1} a_it^i$ in $A \llbrack t \rrbrack$.  Addition and multiplication are given as follows: With $f$ as above and $g(t) = 1 + \sum_{i \geq 1} b_it^i$, their sum and product, respectively, are defined by
\begin{equation*}
\begin{split}
f(t) +_\Lambda g(t) &~=~ 1 + \sum_{i \geq 1}\Biggl[ \sum_{r+s = i}a_r b_s \Biggr] t^i \quad(\text{where } a_0 = b_0 = 1) \\
f(t) \cdot_\Lambda g(t) &~=~ 1 + \sum_{i \geq 1}P_i(a_1, \cdots, a_i; b_1, \cdots, b_i) t^i.
\end{split}
\end{equation*}
Here the $P_i$ are the polynomials which appear in the definition of a $\lambda$-ring (see \S \ref{subsec:filtered lambda-ring}).  If $\alpha \colon A \to B$ is a map of rings, then $\Lambda(\alpha) \colon \Lambda(A) \to \Lambda(B)$ is defined by 
\[
1 + \sum_{i \geq 1}a_i t^i ~\longmapsto~ 1 + \sum_{i \geq 1} \alpha(a_i) t^i.
\]
This makes $\Lambda$ into an endofunctor on the category of rings.

The functor $\Lambda$ is actually (part of) a comonad on the category of rings.  Consider the natural transformations of functors on the category of rings 
\[
\label{eq:eta and lambda}
\eta\colon \Lambda \to \Id \quad \text{and} 
\quad \Lambda_t \colon \Lambda \to \Lambda^2,
\]
which send $f = 1 + \sum_{i \geq 1}a_i t^i \in \Lambda(A)$ to, respectively, $a_1$ and $1 + \sum_{i \geq 1} \lambda^i(f) t^i$.  Here the $\lambda^i$ for $i > 0$ are defined by 
\begin{equation}
\label{eq:lambda operations for Lambda}
\lambda^i(f) ~=~ 1 + \sum_{j \geq 1} P_{j,i}(a_1, \cdots, a_{ij})t^j \in \Lambda(A),
\end{equation}
in which the $P_{j,i}$ are the polynomials appearing in the definition of a $\lambda$-ring.  
The functor $\Lambda$ together with these natural transformations, $\Lambda_t$ and $\eta$, is a comonad on the category of rings.  Moreover, for any ring $A$, $\Lambda(A)$ is a $\lambda$-ring where the $\lambda$-operations are the ones in \eqref{eq:lambda operations for Lambda}.

Now suppose that $A$ is a ring with a certain $\lambda$-ring structure.  Then one can define the following map of rings
\[
\begin{split}
\lambda_t \colon & A ~\to~ \Lambda(A) \\
 & a ~\longmapsto~ \sum_{i \geq 0} \lambda^i(a)t^i
\end{split}
\]
which makes $A$ into a $\Lambda$-coalgebra (that is, $\eta \circ \lambda_t = \Id_A$ and $\Lambda(\lambda_t) \circ \lambda_t = \Lambda_t \circ \lambda_t$).  There is a converse to this.  Suppose that $\xi \colon A \to \Lambda(A)$ is the structure map of a $\Lambda$-coalgebra structure on $A$.  Then there are operations $\lambda^i$ on $A$ defined by the equation
\[
\lambda_t(a) ~=~ \sum_{i \geq 0} \lambda^i(a)t^i ~=~ \xi(a).
\]
The ring $A$ with these operations $\lambda^i$ $(i \geq 0)$ is then a $\lambda$-ring.  Therefore, for any ring $A$, there is a one-to-one correspondence between (resp.\ isomorphism classes of) $\lambda$-ring structures on $A$ and (resp.\ isomorphism classes of) $\Lambda$-coalgebra structures on $A$.   In the rest of this section, this correspondence will be extended to the filtered ring case.

%%%%%%%%%%%%%%%%%%%%%%%%%%%%%%%%%%%%%%%%%%%%%%
%%%%%%%%%%%%%%%%%%%%%%%%%%%%%%%%%%%%%%%%%%%%%%

\subsection{Moduli space of filtered $\lambda$-ring structures via $\Lambda$}
\label{subsec:moduli via Lambda}

In preparation for the main result of this section below, we will first show that the functor $\Lambda$ and the comonad $(\Lambda, \Lambda_t, \eta)$ restrict to the category of filtered rings.

\bigskip
\begin{lemma}
\label{lem1:Lambda}
$\Lambda$ restricts to a functor on the category of filtered rings, where the filtration ideals on $\Lambda(R)$ for a filtered ring $(R, \lbrace I^n \rbrace)$ are
\begin{equation}
\label{eq:filtration for Lambda}
\Lambda(I^n) ~=~ \left\lbrace 1 + \sum_{i \geq 1} a_it^i \in \Lambda(R) \colon a_i \in I^n \text{ for every }i > 0 \right\rbrace.
\end{equation}
\end{lemma}

\begin{proof}
To check that this actually makes $\Lambda(R)$ into a filtered ring, let $f(t) = 1 + \sum_{i \geq 1}a_i t^i$ and $g(t) = 1 + \sum_{i \geq 1}b_i t^i$ be in $\Lambda(I^n)$.  Then 
\[
f +_\Lambda g 
~=~ 1 + \sum_{i \geq 1}c_it^i \quad \text{with} \quad c_i = \sum_{r + s = i}a_rb_s \in I^n.  
\]
(Here we are using the convention $a_0 = b_0 = 1$.)  So $f+_\Lambda g$ lies in $\Lambda(I^n)$.

Now let $h(t) = 1 + \sum_{i \geq 1}d_i t^i$ be in $\Lambda(I^m)$.  Then we have
\[
f\cdot_\Lambda h ~=~ 1 + \sum_{i \geq 1} P_i(a_1,\cdots, a_i; d_1, \cdots, d_i)t^i,
\]
in which each $P_i(a_1, \ldots, a_i;d_1, \ldots, d_i)$ is an integral polynomial with each summand involving at least one of $a_1, \cdots, a_i$ and at least one of $d_1, \cdots, d_i$.  Thus each coefficient $P_i(a_1, \ldots, a_i;d_1, \ldots, d_i)$ lies in $I^{n+m}$, and so $f\cdot_\Lambda h$ lies in $\Lambda(I^{n+m})$.  This shows that $\Lambda(R)$ is a filtered ring with the $\Lambda(I^n)$ as filtration ideals.

Finally, if $\varphi \colon R \to S$ is a filtered ring homomorphism, it is clear that the map $\Lambda(\varphi) \colon \Lambda(R) \to \Lambda(S)$, sending $1 + \sum a_it^i$ to $1 + \sum \varphi(a_i) t^i$, preserves the filtrations on $\Lambda(R)$ and $\Lambda(S)$.

We have shown that $\Lambda$ sends a filtered ring (map) to a filtered ring (map).  Naturality is pretty clear.  So $\Lambda$ restricts to a functor on the category of filtered rings.  This finishes the proof of the lemma.
\end{proof}
\bigskip

From now on whenever we speak of $\Lambda(R)$ as a filtered ring, we mean the filtration in \eqref{eq:filtration for Lambda}.

\bigskip
\begin{thm}
\label{prop1:Lambda}
The natural transformations $\Lambda_t$ and $\eta$ preserve filtrations when applied to filtered rings.  Therefore, $(\Lambda, \Lambda_t, \eta)$ restricts to a comonad on the category of filtered rings.
\end{thm}

\begin{proof}
In view of Lemma \ref{lem1:Lambda} and the fact that $(\Lambda, \Lambda_t, \eta)$ is already a comonad on the category of rings, it suffices to prove the first assertion, which is pretty clear for $\eta$.  To see that $\Lambda_t$ is a filtered ring map when applied to filtered rings $(R, \lbrace I^n \rbrace)$, one only has to note that the $P_{j,i}$ are integral polynomials without constant terms, and so $\lambda^i(f)$ $(i > 0)$ lies in $\Lambda(I^n)$ whenever $f \in \Lambda(R)$ does.
\end{proof}

This proof also yields the following result, which is a filtered analog of the fact that, for any ring $R$, $\Lambda(R)$ has a canonical $\lambda$-ring structure.

\bigskip
\begin{cor}
\label{cor2:Lambda}
Let $(R, \lbrace I^n \rbrace)$ be a filtered ring.  Then the filtered ring $\Lambda(R)$ with the operations $\lambda^i$ as in \eqref{eq:lambda operations for Lambda} is a filtered $\lambda$-ring.
\end{cor}

\begin{proof}
Indeed, what remains to be proved is that $\lambda^i(\Lambda(I^n)) \subseteq \Lambda(I^n)$ for any $i > 0$, $n > 0$, which was already observed in the proof of Theorem \ref{prop1:Lambda}.
\end{proof}

The next result shows that the filtered ring $\Lambda(R)$ has topological properties similar to those of $R$; namely, one of them is complete (resp.\ Hausdorff) if and only if the other is so.  Recall that a filtered ring $(R, \lbrace I^n \rbrace)$ is said to be complete (resp.\ Hausdorff) if the natural map 
\[
p \colon R ~\to~ \varprojlim R/I^n
\]
of filtered rings is surjective (resp.\ injective).  Note that this map $p$ is injective if and only if the intersection $\cap_{n \geq 1} I^n$ is $0$.  This has relevance, for example, in topology: If $X$ is a space in the genus of a torsionfree classifying space $BG$ of a simply-connected compact Lie group $G$, then its integral $K$-theory $K(X)$ is complete Hausdorff (see \cite{yau1}).

\bigskip
\begin{prop}
\label{prop:complete and hausdorff}
Let $(R, \lbrace I^n \rbrace)$ be a filtered ring.  Then $(R, \lbrace I^n \rbrace)$ is complete $($resp.\ Hausdorff$)$ if and only if $(\Lambda(R), \lbrace \Lambda(I^n) \rbrace)$ is complete $($resp.\ Hausdorff$)$.
\end{prop}

\begin{proof}
From the definition of $\Lambda(I^n)$, it is easy to see that the condition 
\[
\bigcap_{n \geq 1} I^n ~=~ (0) ~\subset~ R 
\]
is equivalent to the condition
\[
\bigcap_{n \geq 1} \Lambda(I^n) ~=~ (0) ~\subset~ \Lambda(R).  
\]
That is, $R$ is Hausdorff if and only if $\Lambda(R)$ is Hausdorff.

For the assertion about completeness, consider the commutative diagram of filtered rings and filtered ring maps:
\begin{equation}
\xymatrix{ \Lambda(R) \ar[r]^-\gamma \ar[dr]_\alpha & \varprojlim \Lambda(R)/\Lambda(I^n) \ar[d]_\beta \\
& \varprojlim \Lambda(R/I^n)}
\end{equation}
%\[
%\begin{diagram}
% \node{\Lambda(R)} \arrow{e,t}{\gamma} \arrow{se,b}{\alpha} 
% \node{\varprojlim \Lambda(R)/\Lambda(I^n)} \arrow{s,t}{\beta} \\
% \node{} \node{\varprojlim \Lambda(R/I^n)}
%\end{diagram}
%\]
The map $\gamma$ is induced by the projection maps $\Lambda(R) \to \Lambda(R)/\Lambda(I^n)$ and $\alpha$ and $\beta$ by $R \to R/I^n$.  It is clear that for each $n$, the map $\Lambda(R) \to \Lambda(R/I^n)$ is surjective with kernel $\Lambda(I^n)$; thus $\beta$ is an isomorphism.  Moreover, it is evident that $p \colon R \to \varprojlim R/I^n$ is surjective if and only if $\alpha$ is so.  The fact that $\beta$ is an isomorphism now implies that $p$ is surjective if and only if $\gamma$ is so.  That is, $R$ is complete if and only if $\Lambda(R)$ is so.

This finishes the proof of the proposition.
\end{proof}

Now we are ready for the main result of this section, which is that the correspondence between $\lambda$-ring structures and $\Lambda$-coalgebra structures over a ring has a filtered analog.

\bigskip
\begin{thm}
\label{thm:Lambda}
Let $(R, \lbrace I^n \rbrace)$ be a filtered ring.  Then the set of filtered $\lambda$-ring structures on $R$ is canonically isomorphic to the set of $\Lambda$-coalgebra structures on $R$, where $\Lambda$ is regarded as a comonad on the category of filtered rings.
\end{thm}

\begin{proof}
First suppose that $R$ is given a filtered $\lambda$-ring structure.  Then for an element $r$ in $I^n$, $\lambda^i(r)$ belongs to $I^n$ for every $i \geq 1$, so that 
\[
\lambda_t(r) ~=~ 1 + \sum_{i \geq 1} \lambda^i(r)t^i 
\]
lies in $\Lambda(I^n)$.  Therefore, the map $\lambda_t \colon R \to \Lambda(R)$ is a filtered ring homomorphism.  Since this map (between rings) already makes $R$ into a $\Lambda$-coalgebra when $\Lambda$ is regarded as a comonad on the category of rings, it also gives rise to a $\Lambda$-coalgebra structure on the filtered ring $R$ when $\Lambda$ is regarded as a comonad on the category of filtered rings.

Conversely, suppose that $\xi \colon R \to \Lambda(R)$ is the structure map of a $\Lambda$-coalgebra structure on the filtered ring $R$.  Then it yields a $\lambda$-ring structure on $R$ by declaring that $\xi = \lambda_t$.  To see that it is actually a filtered $\lambda$-ring, we must show that $\lambda^i(I^n) \subseteq I^n$ for any $i > 0$, $n > 0$.  Now if $r$ is in $I^n$, then 
\[
\xi(r) ~=~ \lambda_t(r) ~=~ 1 + \sum_{i \geq 1} \lambda^i(r) t^i 
\]
lies in $\Lambda(I^n)$.  Thus, $\lambda^i(r)$ belongs to $I^n$ for every $i \geq 1$.  This shows that the $\lambda$-ring structure on $R$ given by $\xi$ is a filtered $\lambda$-ring.

This finishes the proof of the theorem.
\end{proof}

By slightly modifying the above proof, one can also classify the isomorphism classes of filtered $\lambda$-ring structures over a filtered ring.

\bigskip
\begin{cor}
\label{cor1:Lambda}
The set of isomorphism classes of filtered $\lambda$-ring structures on a filtered ring $R$ is canonically isomorphic to the set of isomorphism classes of $\Lambda$-coalgebra structures on $R$, where $\Lambda$ is regarded as a comonad on the category of filtered rings.
\end{cor}

If $R$ is a filtered ring, then the filtration generates a topology on $R$, so in this case both $R$ and $\Lambda(R)$ are topological spaces.  The set of continuous maps from $R$ to $\Lambda(R)$ is then a topological space with the compact-open topology.  Now Theorem \ref{thm:Lambda} yields the following result, which is what the title of this paper refers to.

\bigskip
\begin{cor}
\label{cor1:moduli space}
Regard $\Lambda$ as a comonad on the category of filtered rings.  Then the set of filtered $\lambda$-ring structures on a filtered ring $R$ has a canonical topology as the subspace of the space of continuous maps from $R$ to $\Lambda(R)$ consisting of the filtered ring homomorphisms which make $R$ into a $\Lambda$-coalgebra.
\end{cor}

For a filtered ring $R$, the set of filtered $\lambda$-ring structures on $R$ with the topology in Corollary \ref{cor1:moduli space} is called the {\it moduli space of filtered} $\lambda$-{\it ring structures on} $R$.

%%%%%%%%%%%%%%%%%%%%%%%%%%%%%%%%%%%%%%%%%%%%%%%%%%%%%%%%%%
%%%%%%%%%%%%%%%%%%%%%%%%%%%%%%%%%%%%%%%%%%%%%%%%%%%%%%%%%%

\subsection{Moduli space of filtered $\lambda$-ring structures via $\bW$}
\label{subsec:Witt}

We now show that the results in the previous section can be restated with the functor $\bW$ in place of $\Lambda$.  More information about the functor $\bW$ (in the category of rings) can be found in Hazewinkel's book \cite[Ch.\ 3]{hazewinkel}.

For a ring $A$, the universal ring of Witt vectors on $A$, denoted $\bW(A)$, is the ring defined as follows.  As a set $\bW(A)$ is the countable product $\prod_{i=1}^\infty A$.  There are functions 
\[
w_n \colon \bW(A) ~\to~ A \quad (n \geq 1)
\]
defined by, for an element $\mathbf{a} = (a_i)$ in $\bW(A)$, the formula
\[
w_n(\mathbf{a}) ~=~ \sum_{d \vert n}\, d a_d^{\frac{n}{d}}.
\]
If $\mathbf{a} = (a_i)$ and $\mathbf{b} = (b_i)$ are elements of $\bW(A)$, then their sum and product in the ring $\bW(A)$ are the unique elements $\mathbf{c} = \mathbf{a} +_W \mathbf{b}$ and $\mathbf{d} = \mathbf{a}\cdot_W \mathbf{b}$, respectively, such that for any $n \geq 1$ the following equations hold:
\begin{equation*}
\begin{split}
w_n(\mathbf{c}) &~=~ w_n(\mathbf{a}) + w_n(\mathbf{b}) \\
w_n(\mathbf{d}) &~=~ w_n(\mathbf{a}) w_n(\mathbf{b}).
\end{split}
\end{equation*}
If $\varphi \colon A \to B$ is a ring map, then $\bW(f) \colon \bW(A) \to \bW(B)$ is defined by sending $(a_i)$ to $(\varphi(a_i))$.  Then $\bW$ is an endofunctor and the $w_n \colon \bW \to \Id$ are natural transformations of functors on the category of rings.

The functor $\bW$ is related to $\Lambda$ via the following map, called the exponential isomorphism:
\[
\begin{split}
E(A) \colon \bW(A) & ~\xrightarrow{\cong}~ \Lambda(A) \\
(a_1, a_2, \cdots) & ~\mapsto~ \prod_{i \geq 1}(1 + a_it^i)
\end{split}
\]
This is a ring isomorphism for any ring $A$.  Actually $E \colon \bW \to   \Lambda$ is a natural transformation of functors on the category of rings.  This natural transformation implies that $\bW$ is also a comonad on the category of rings, in which the counit is $\eta E$ and the comultiplication is $(E^2)^{-1} \Lambda_t E$.  Here $E^2 \colon W^2 \to \Lambda^2$ is the composition
\[
W^2 ~\xrightarrow[\cong]{E}~ \Lambda W ~\xrightarrow[\cong]{\Lambda E}~ \Lambda^2.
\]
The natural transformation $E$, in fact, gives an isomorphism 
\[
(\Lambda, \Lambda_t, \eta) ~\xrightarrow{\cong}~ (\bW, (E^2)^{-1} \Lambda_t E, \eta E)
\]
of comonads on the category of rings.

The following result will imply immediately that the results in the previous section, in particular Theorem \ref{thm:Lambda} and its corollaries, are still valid with $\Lambda$ replaced by $\bW$ everywhere.  
For a filtered ring $(R, \lbrace I^n \rbrace)$, denote by $\bW(I^n)$ the subset of $\bW(R)$ consisting of elements $(r_i)$ with each $r_i \in I^n$.

\bigskip
\begin{thm}
\label{lem2:Lambda}
Let $(R, \lbrace I^n \rbrace)$ be a filtered ring and let ${\bf r} = (r_i)$ be an element of $\bW(R)$.  Then ${\bf r}$ lies in $\bW(I^n)$ if and only if its image under $E(R)$ lies in $\Lambda(I^n)$.  Therefore, $\bW(R)$ with the ideals $\bW(I^n)$ is a filtered ring and $\bW$ and $(\bW, (E^2)^{-1} \Lambda_t E, \eta E)$ restrict to, respectively, a functor and a comonad on the category of filtered rings.  Moreover, $E$ is an isomorphism of comonads on the category of filtered rings and $E(R)$ is a filtered ring isomorphism.
\end{thm}

\begin{proof}
Note that the image of ${\bf r}$ under $E(R)$ is
\[
\label{eq1:lem2}
E(R)({\bf r}) ~=~ \prod_{i \geq 1}(1 + r_it^i) 
~=~ 1 + r_1t + r_2t^2 + Q_3t^3 + Q_4t^4 + \cdots.
\]
Here for $n \geq 3$ the coefficient $Q_n$ has the form:
\begin{equation}
\label{eq2:lem2}
Q_n ~=~ \left(\sum r_{i_1}\cdots r_{i_k}\right) + r_n,
\end{equation}
where the sum inside the parentheses is taken over all partitions, $i_1 + \cdots + i_k = n$, of $n$ with $1 \leq i_1 < i_2 < \cdots < i_k < n$.  This shows that if ${\bf r}$ lies in $\bW(I^n)$, then its image under $E(R)$ lies in $\Lambda(I^n)$.

Conversely, suppose that $E(R)({\bf r})$ lies in $\Lambda(I^n)$.  From the definition of the exponential isomorphism $E$, it is easy to see that the elements $r_1$ and $r_2$ belong to $I^n$.  Suppose by induction that $r_j$ lies in $I^n$ for every $j < m$.  Then it follows from \eqref{eq2:lem2} that $Q_m - r_m$ lies in $I^n$, which implies that $r_m$ lies in $I^n$ since $Q_m$ does.

This finishes the proof of the first assertion.  The other assertions follow from this, Lemma \ref{lem1:Lambda}, and the two paragraphs immediately preceding this theorem. 
\end{proof}

                                   %%  Power series  %% 

\section{Filtered $\lambda$-rings over power series rings}
\label{sec:power series}

The purpose of this section is to study (co)representability of the functor which sends a ring $R$ to the set of filtered $\lambda$-ring structures over $R\llbrack x \rrbrack$.  The main result is Theorem \ref{thm:co-representing ring}, which implies Theorem \ref{thm-intro:Lazard} in the Introduction.  Throughout this section, the indeterminate $x$ is given a fixed positive filtration $d$, and $R \llbrack x \rrbrack$ is given the $x$-adic filtration.  First we recall a result of Wilkerson about Adams operations determining a unique $\lambda$-ring structure, which we will need later on in this section and also in the next section.

%%%%%%%%%%%%%%%%%%%%%%%%%%%%%%%%%%
%%%%%%%%%%%%%%%%%%%%%%%%%%%%%%%%%%

\subsection{Wilkerson's Theorem}
\label{subsec:wilkerson's theorem}

Recall that a $\lambda$-ring $R$ has Adams operations $\psi^n \colon R \to R$ ($n \geq 1$) defined by the Newton formula
\begin{equation}
\label{eq:Newton}
\psi^n(r) - \lambda^1(r)\psi^{n-1}(r) + \cdots + (-1)^{n-1}\lambda^{n-1}(r) \psi^1(r) + (-1)^nn\lambda^n(r) = 0.
\end{equation} 
They satisfy the following properties:
\begin{itemize}
\item all the $\psi^n$ are $\lambda$-ring maps,
\item $\psi^1 = \Id$,
\item $\psi^m \psi^n = \psi^{mn} = \psi^n \psi^m$,
\item $\psi^p(r) \equiv r^p$ (mod $pR$) for each prime $p$ and element $r$ in $R$.
\end{itemize}
If $R$ is a filtered $\lambda$-ring, then the corresponding Adams operations are filtered $\lambda$-ring maps.  There is a partial converse to this.  Suppose that $(R, \lbrace I^n \rbrace)$ is a filtered ring with a $\lambda$-ring structure in which the Adams operations are filtered ring maps.  If the filtration ideals $I^m$ are all \emph{closed under dividing by integers} (that is, $nr \in I^m$ implies $r \in I^m$ for $n, m > 0$, $r \in R$), then the Newton formula \ref{eq:Newton} implies that they are also closed under the operations $\lambda^i$ for $i > 0$.  For example, the power series ring $R \llbrack x \rrbrack$ with the $x$-adic filtration has the property that each filtration ideal is closed under dividing by integers.

Adams operations are easier to work with because they are (filtered) ring maps.  So it is tempting to ask if a ring $R$ equipped with endomorphisms $\psi^n$ ($n \geq 1$) satisfying the last three properties above actually has a $\lambda$-ring structure with the $\psi^n$ as Adams operations.  The answer is yes, as the following theorem of Wilkerson \cite[Prop.\ 1.2]{wilkerson} shows, provided that $R$ is torsionfree as a $\bZ$-module.

\bigskip
\begin{thm}[Wilkerson]
\label{thm:wilkerson}
Suppose that $R$ is a torsionfree ring equipped with ring homomorphisms $\psi^n \colon R \to R$ for $n \geq 1$ satisfying the properties: 
   \begin{enumerate}
   \item $\psi^1 = \Id$ and $\psi^m\psi^n = \psi^{mn}$,
   \item $\psi^p(r) \equiv r^p$ $($mod $pR)$ for each prime $p$ and element $r$ in $R$.
   \end{enumerate}
Then there is a unique $\lambda$-ring structure over $R$ with the given $\psi^n$ as Adams operations.
\end{thm}

Since in this section our attention is on Adams operations rather than $\lambda$-operations, the following notion will be useful.

%%%%%%%%%%%%%%%%%%%%%%%%%%
%%%%%%%%%%%%%%%%%%%%%%%%%%

\subsection{Filtered $\psi$-ring}
\label{subsec:filtered psi-ring}

A $\psi$-\emph{ring} is a ring $R$ equipped with ring homomorphisms $\psi^p \colon R \to R$ ($\psi^p_R$ if we wish to specify $R$), $p$ primes, such that 
\begin{itemize}
\item $\psi^p(r) \equiv r^p$ (mod $pR$) for each element $r$ in $R$ and prime $p$, 
\item $\psi^p\psi^q = \psi^q \psi^p$ for any primes $p$ and $q$.
\end{itemize}
A $\psi$-ring map $f \colon R \to S$ between two $\psi$-rings is a ring homomorphism $f \colon R \to S$ which commutes with the structure maps: $f \psi^p_R = \psi^p_S f$ for each prime $p$.

A \emph{filtered} $\psi$-\emph{ring} is a $\psi$-ring $R$ which is also a filtered ring for which the structure maps $\psi^p$ are filtered ring homomorphisms.  A map between two filtered $\psi$-rings is a filtered ring homomorphism which is also a $\psi$-ring map.

Any filtered $\lambda$-ring gives rise to a filtered $\psi$-ring by using its Adams operations $\psi^p$.  Conversely, Wilkerson's Theorem \ref{thm:wilkerson} implies that a filtered $\psi$-ring structure on a torsionfree filtered ring $R$ in which the filtration ideals are closed under dividing by integers, determines a unique filtered $\lambda$-ring structure on $R$.

The following preliminary result gives necessary and sufficient conditions in order that $R \llbrack x \rrbrack$ be a filtered $\psi$-ring.

\bigskip
\begin{lemma}
\label{lem1:filtered ring}
Let $R$ be a ring between $\bZ$ and $\bQ$.  Then a filtered $\psi$-ring $(=$ filtered $\lambda$-ring$)$ structure over $R\llbrack x \rrbrack$ determines and is determined by power series $\psi^p(x)$ in $R\llbrack x \rrbrack$, one for each prime $p$, such that for any primes $p$ and $q$ the following conditions hold:
\begin{itemize}
\item $\psi^p(0) = 0$,
\item $\psi^p(x) \equiv x^p$ $($mod $pR)$,
\item $\psi^p (\psi^q(x)) = \psi^q (\psi^p(x))$.
\end{itemize}
\end{lemma}

\begin{proof}
The structure maps of a filtered $\psi$-ring structure on $R \llbrack x \rrbrack$ clearly gives rise to such a sequence of power series.  

Conversely, suppose given a sequence of power series $\psi^p(x)$ satisfying the above properties.  Since any ring endomorphism on $R \llbrack x \rrbrack$ must restrict to the identity map on $R$, the power series  $\psi^p(x)$ can be extended uniquely to filtered ring endomorphisms $\psi^p$ on $R \llbrack x \rrbrack$.  It is clear that these filtered ring maps satisfy 
\[
\psi^p \psi^q ~=~ \psi^q \psi^p
\]
for any primes $p$ and $q$.   Also, the hypothesis about $\psi^p(x)$ implies that 
\[\psi^p(\alpha) ~\equiv~ \alpha^p \qquad \left(\text{mod } pR\llbrack x \rrbrack \right)
\]
for any power series $\alpha$ over $R$, since this is true when $\alpha$ is in $R$ and $\psi^p$ is the identity map on $R$.  Therefore, the power series $\psi^p(x)$ ($p$ prime) determine a unique filtered $\psi$-ring structure on $R \llbrack x \rrbrack$.
\end{proof}

%%%%%%%%%%%%%%%%%%%%%%%%%%%%%%%%
%%%%%%%%%%%%%%%%%%%%%%%%%%%%%%%%

\subsection{The representability question}
\label{subsec:representability}

Before we state the main result of this section, note that if $f \colon R \to S$ is a ring homomorphism and if $R\llbrack x \rrbrack$ is equipped with a filtered $\psi$-ring structure in which each of the structure maps $\psi^p_R$ restricts to the identity on $R$, then it can be pushed forward to a filtered $\psi$-ring structure $\lbrace \psi^p_S \rbrace$ over $S \llbrack x \rrbrack$ as follows.  (Note that we are using the subscripts $R$ and $S$ for filtered $\psi$-ring structures over $R \llbrack x \rrbrack$ and $S \llbrack x \rrbrack$, respectively.)  Denote by $f_*(\psi^p_R(x))$ the power series over $S$ obtained by applying $f$ to each coefficient of $\psi^p_R(x)$, and take $\psi^p_S$ to be the unique filtered ring endomorphism on $S \llbrack x \rrbrack$ which restricts to the identity map on $S$ and satisfies
\[
\psi^p_S(x) ~=~ f_*(\psi^p_R(x)).
\]
In this case we also write $f_*\left(\psi^p_R\right)$ for the map $\psi^p_S$.

For two rings $R$ and $S$, the set of ring homomorphisms from $R$ to $S$ is denoted by $\Ring(R,S)$.

We are now ready for the main result of this section, which implies Theorem \ref{thm-intro:Lazard} in the Introduction.

\bigskip
\begin{thm}
\label{thm:co-representing ring}
There exist a binomial domain $U$ and a filtered $\psi$-ring structure $\left\lbrace\psi^p_{univ}\right\rbrace$ over $U \llbrack x \rrbrack$ such that the following statements hold.
\begin{enumerate}
\item Each $\psi^p_{univ}$ restricts to the identity map on $U$.
\item For every ring $R$ between $\bZ$ and $\bQ$, the natural map 
\[
\begin{split}
\Psi_R \colon  \Ring(U,R&)  ~\xrightarrow~ \left\lbrace\text{filtered }\lambda \text{-ring structures over } R\llbrack x \rrbrack \right\rbrace \\
 & f \longmapsto \lbrace f_*(\psi^p_{univ})\rbrace
\end{split}
\]
is an isomorphism.
\end{enumerate}
\end{thm}

\begin{proof}
The ring $U$ is defined to be the quotient
\[
\label{eq:the ring U}
U ~=~ \frac{\bZ \left\lbrack \lbrace v_{(p,\, i, q_1,\, \ldots,\, q_n)} \rbrace \right\rbrack /J}{(\bZ\text{-torsions})}.
\]
In other words, take the polynomial ring $\bZ \left\lbrack \lbrace v_{(p,\, i, q_1,\, \ldots,\, q_n)} \rbrace \right\rbrack$, divide out by the ideal $J$ (to be defined below), and then divide out by the ideal generated by the $\bZ$-torsions in the resulting quotient $\bZ \left\lbrack \lbrace v_{(p,\, i, q_1,\, \ldots,\, q_n)} \rbrace \right\rbrack /J$.  Here the $v_{(p,\, i, q_1,\, \ldots,\, q_n)}$ are independent indeterminates, where $n \geq 0$, $i \geq 1$, and $p$, $q_1$, $\ldots$ , $q_n$ any primes.  The ideal $J$ is generated by the elements $w_{(p,\, q,\, l)}$ and $V(p, i, q_1, \ldots, q_{n+1})$ defined as follows.  We first define
\[
\label{eq:definition of the u_{(p,i)}}
  u_{(p,\, i)} ~=~ \begin{cases} pv_{(p,\, i)} & \text{if } i \not=p\\
    1 + pv_{(p,\, p)} & \text{if } i = p.\end{cases}
\]
The elements $w_{(p, \, q,\, l)}$ are then defined by the equation
\begin{equation}
\label{eq:commutativity}
\begin{split}
\sum_{i \geq 1}\, u_{(p,\, i)}\Biggl[\sum_{j \geq 1}\, u_{(q,\, j)}x^j\Biggr]^i ~-~
\sum_{i \geq 1}\, u_{(q,\, i)}\Biggl[\sum_{j \geq 1}\, u_{(p,\, j)}x^j\Biggr]^i ~=~ \sum_{l \geq 1}\, w_{(p,\, q,\, l)} x^l.
\end{split}
\end{equation}
The elements $V(p, i, q_1, \ldots, q_{n+1})$ are defined to be
\begin{equation}
\label{eq:the elements V}
  \begin{split}
V(p, i, q_1, \ldots, q_{n+1}) ~=~ & v_{(p,\, i, q_1,\, \ldots,\, q_n)}^{q_{n+1}} ~-~ v_{(p,\, i, q_1,\, \ldots,\, q_n)} \\
 & ~-~ q_{n+1}v_{(p,\, i, q_1,\, \ldots,\, q_n,\, q_{n+1})}.
  \end{split}
\end{equation}
Here $p$, $q$, $q_1$, $\ldots$ , $q_{n+1}$ are primes, $n \geq 0$, and $l, i \geq 1$.

It follows immediately from the definition of the elements $V(p, i, q_1, \ldots, q_{n+1})$ that each element $u$ in $U$ satisfies the congruence condition 
\begin{equation}
\label{eq:congruence}
u^p ~\equiv~ u \quad (\text{mod }pU).
\end{equation}
Since $U$ is also torsionfree, Wilkerson's Theorem \ref{thm:wilkerson} implies that it is a binomial domain.

Now consider the following power series over $U \llbrack x \rrbrack$:
\[
\label{eq:adams operations for U}
\psi^p_{univ}(x) 
~=~  \sum_{i \geq 1}\, u_{(p,i)}x^i
~=~ (1 + pv_{(p,\, p)})x^p 
~+~ 
\sum_{i \geq 1, \, i \not= p}\, pv_{(p,\, i)}x^i. 
\]
These power series $\psi^p_{univ}(x)$ extend uniquely to filtered ring endomorphisms on $U \llbrack x \rrbrack$ if one insists that each $\psi^p_{univ}$ restricts to the identity on $U$.  To see that $\left\lbrace \psi^p_{univ} \right\rbrace$ is a filtered $\psi$-ring structure on $U \llbrack x \rrbrack$, first note that \eqref{eq:commutativity} implies that
\[
\psi^p_{univ}\left(\psi^q_{univ}(x)\right) ~=~
\psi^q_{univ}\left(\psi^p_{univ}(x)\right)
\]
for any primes $p$ and $q$.  So $\psi^p_{univ}\psi^q_{univ}$ and $\psi^q_{univ}\psi^p_{univ}$ coincide as filtered ring endomorphisms.  Moreover, since one has that
\[
\psi^p_{univ}(x) ~\equiv~ x^p \quad (\text{mod }pU\llbrack x \rrbrack),
\]
the condition
\[
\psi^p_{univ}(\alpha) ~\equiv~ \alpha^p \quad (\text{mod }pU\llbrack x \rrbrack) \quad \text{for every }\alpha \in U\llbrack x \rrbrack
\]
is equivalent to \eqref{eq:congruence}, which we already know is true.  Therefore, $\lbrace \psi^p_{univ} \rbrace$ is a filtered $\psi$-ring structure on $U\llbrack x \rrbrack$ in which each structure map restricts to the identity map on $U$.  This proves the first assertion.

Now we consider the second assertion.  To prove the surjectivity of $\Psi_R$, let $R$ be a ring between $\bZ$ and $\bQ$ with a given filtered $\psi$-ring structure $\lbrace \psi^p \rbrace$ on $R \llbrack x \rrbrack$ which is determined by the power series (see Lemma \ref{lem1:filtered ring}) 
\[
\psi^p(x) ~=~ \sum_{i \geq 1}\, r_{(p,\, i)}x^i.  
\]
We want a ring map $f \colon U \to R$ which carries the filtered $\psi$-ring structure $\left\lbrace \psi^p_{univ} \right\rbrace$ on $U \llbrack x \rrbrack$ to the given one on $R \llbrack x \rrbrack$.  We first specify the values of $f$ on the generators $v_{(p,\, i)}$ to be
\[
f(v_{(p,\, i)}) ~=~ 
\begin{cases} r_{(p,\, i)}/p & \text{ if } i \not= p \\
\left(r_{(p,\, p)} - 1\right)/p & \text{ if } i = p. 
\end{cases}
\] 
This makes sense because 
\[
r_{(p,\, i)} ~\equiv~
\begin{cases}0 & (\text{mod } pR) \quad \text{if } i \not= p \\
  1 & (\text{mod }pR) \quad \text{if } i = p. 
\end{cases}
\]
The values of $f$ on the other generators of $U$ are then inductively determined by the formula:
\begin{equation}
\label{eq:value of f}
f\left(v_{(p,\, i,\, q_1,\, \ldots,\, q_{n+1})}\right) ~=~
\frac{f\left(v_{(p,\, i,\, q_1,\, \ldots,\, q_n)}\right)^{q_{n+1}} - f\left(v_{(p,\, i,\, q_1,\, \ldots,\, q_n)}\right)}{{q_{n+1}}}.
\end{equation}
This shows that the map $\Psi_R$ is surjective.

To show that $\Psi_R$ is injective, suppose that $f$ and $g$ are two ring maps from $U$ to $R$ such that \[
f_*(\psi^p_{univ}(x)) ~=~ g_*(\psi^p_{univ}(x))
\]
for each prime $p$.  We must show that $f = g$.  The hypothesis implies that 
\[
f(u_{(p,\, i)}) ~=~ g(u_{(p,\, i)})
\]
for each prime $p$ and $i \geq 1$.  Since $R$ is torsionfree, this in turn implies that 
\[
f(v_{(p,\, i)}) ~=~ g(v_{(p,\, i)}).  
\]
The formula \eqref{eq:value of f} now shows that $f$ and $g$ must agree on all the generators $v_{(p, i, q_1, \ldots, q_n)}$ of $U$, whence $f = g$.

This finishes the proof of the theorem.
\end{proof}

                                    %%  Examples  %%

\section{Examples}
\label{sec:examples}

In this final section we study points in the moduli spaces of filtered $\lambda$-ring structures over certain filtered rings.  In each case we will exhibit uncountably many mutually non-isomorphic filtered $\lambda$-ring structures, which in particular implies that the moduli space has uncountably many points.  Unlike in \cite{yau1} where the existence of uncountably many mutually non-isomorphic $\lambda$-ring structures over the power series ring $\bZ \llbrack x_1, \ldots, x_n\rrbrack$ was obtained by a combination of \emph{topological} results, our method here is purely algebraic.  We will write down directly Adams operations and then use Wilkerson's Theorem \ref{thm:wilkerson} to lift them to a $\lambda$-ring structure.

We begin with

%%%%%%%%%%%%%%%%%%%%%%%%%%%%%%%%
%%%%%%%%%%%%%%%%%%%%%%%%%%%%%%%%

\subsection{Rings of dual numbers}
\label{subsec:dual numbers}

Recall that a binomial domain is a torsionfree ring $B$ for which the binomial symbol
\[
\binom{b}{n} = \frac{b(b-1) \cdots (b-n+1)}{n!} \in B \otimes \bQ
\]
lies in $B$ for any element $b$ in $B$ and positive integer $n$.  Examples of binomial domains include the integers, the $p$-local integers, the $p$-adic integers, and $\bZ$-torsionfree $\bQ$-algebras (and hence fields of characteristic $0$).

Let $\varepsilon$ be a square-zero variable: $\varepsilon^2 = 0$.  The ring $B \lbrack \varepsilon \rbrack$, known as the ring of dual numbers over $B$, is given the $\varepsilon$-adic filtration; that is, $I^0 = B \lbrack \varepsilon \rbrack$, $I^1 = (\varepsilon)$, and $I^n = (0)$ for any $n \geq 2$.  Note that $B \lbrack \varepsilon \rbrack$, as a $B$-module, is free of rank $2$, and in particular it can be a finitely generated abelian group.

We now show that, despite the ``small'' size of $B \lbrack \varepsilon \rbrack$ as a $B$-module, its moduli space of filtered $\lambda$-ring structures has many points.

\bigskip
\begin{thm}
\label{thm:dual numbers}
The filtered ring $B \lbrack \varepsilon \rbrack$ of dual numbers over any non-zero binomial domain $B$ admits uncountably many mutually non-isomorphic filtered $\lambda$-ring structures.
\end{thm}

\begin{proof}
Let $R$ denote the filtered ring $B \lbrack \varepsilon \rbrack$.  We will first show that any sequence $(a_p)$ of elements in $B$ indexed by the set of primes with each $a_p$ $p$-divisible in $B$ gives rise to a filtered $\lambda$-ring structure on $R$.  So suppose that $(a_p)$ is such a sequence.  Define an associated sequence of filtered-ring endomorphisms, $\psi^n \colon R \to R$ $(n \geq 1)$, as follows.  Set $\psi^1 = \Id$ and if $n$ has prime decomposition $n = p_1^{e_1} \cdots p_k^{e_k}$, then for any elements $\alpha, \beta \in B$, we set
\[
\psi^n(\alpha + \beta \varepsilon) 
~=~ \alpha + \beta a_{p_1}^{e_1} \cdots a_{p_k}^{e_k} \varepsilon.
\] 
It is clear that the $\psi^n$ are filtered-ring endomorphisms such that 
\[
\psi^n \psi^m ~=~ \psi^m \psi^n ~=~ \psi^{mn}.
\]
Also, for any prime $p$ we have that 
\[
\begin{split}
\psi^p(\alpha + \beta \varepsilon) 
&~=~ \alpha + \beta a_p \varepsilon  \\
&~\equiv~ \alpha \hspace{.65in} (\text{mod } pR) \\
&~\equiv~ \alpha^p + \beta^p \varepsilon^p \quad (\text{mod } pR) \\
&~\equiv~ (\alpha + \beta \varepsilon)^p \quad (\text{mod } pR)
\end{split}
\]
It then follows easily from Wilkerson's Theorem \ref{thm:wilkerson} that there is a unique filtered $\lambda$-ring structure on $R$ with these $\psi^n$ as Adams operations.

To finish the proof, we will now show that filtered $\lambda$-ring structures on $R$ arising from different sequences are not isomorphic.  So 
suppose that two filtered $\lambda$-ring structures $S_1$ and $S_2$ over $R$ corresponding to sequences $(a_p)$ and $(b_p)$, respectively, are isomorphic as filtered $\lambda$-rings.  We claim that $a_p = b_p$ for each $p$.  Indeed, if $\varphi \colon S_1 \to S_2$ is a filtered $\lambda$-ring isomorphism, then $\varphi$ must send $\varepsilon$ to $u\varepsilon$ for some unit $u \in B^*$.  Therefore, since \[
\varphi\psi^p ~=~ \psi^p \varphi,
\]
 we have that 
\[
ua_p\varepsilon ~=~ ub_p\varepsilon,
\] 
whence $a_p = b_p$.

Now since $B$ is a non-zero, characteristic $0$ ring, there are uncountably many distinct sequences $(a_p)$ in it with each $a_p$ $p$-divisible.  Thus, together with what we just proved above, we conclude that there are uncountably many distinct isomorphism classes of filtered $\lambda$-ring structures on $R$.

This finishes the proof of the theorem.
\end{proof}

Note that the above argument actually gives a complete classification of all of the isomorphism classes of filtered $\lambda$-ring structures when $B$ is any ring between $\bZ$ and $\bQ$.  To see this, notice that if $\psi^p$ is the Adams operations of a filtered $\lambda$-ring structure over $B\lbrack \varepsilon \rbrack$, then 
\[
\psi^p(m/n) ~=~ m/n 
\]
for any element $m/n \in B$.  Now since we have that 
\[
\psi^p(\varepsilon) ~=~ a_p\varepsilon ~\equiv~ \varepsilon^p ~=~ 0 \quad (\text{mod } pB) 
\]
for some element $a_p$ in $B$, so $a_p$ must be $p$-divisible in $B$.  This means that for any elements $\alpha$ and $\beta$ in $B$, one has that 
\[
\psi^p(\alpha + \beta \varepsilon) 
~=~ \psi^p(\alpha) + \psi^p(\beta \varepsilon) 
~=~  \alpha + \beta a_p \varepsilon.
\]
Now if $n$ has prime decomposition, $n = p_1^{e_1} \cdots p_k^{e_k}$, then the fact that $\psi^r \psi^s = \psi^s \psi^r = \psi^{rs}$ for any $r$ and $s$ implies that
\[
\psi^n(\alpha + \beta \varepsilon) 
~=~ (\psi^{p_1})^{e_1} \cdots (\psi^{p_k})^{e_k}(\alpha + \beta \varepsilon) 
~=~ \alpha + \beta a_{p_1}^{e_1} \cdots a_{p_k}^{e_k} \varepsilon.
\] 
Thus, the given filtered $\lambda$-ring structure over $B \lbrack \varepsilon \rbrack$ is among those constructed in the proof of Theorem \ref{thm:dual numbers}.  Let us state this as a separate corollary.

\bigskip
\begin{cor}
\label{cor:dual numbers}
Let $B$ be any ring between $\bZ$ and $\bQ$.  Then a filtered $\lambda$-ring structure on the filtered ring $B \lbrack \varepsilon \rbrack$ of dual numbers must have, for any prime $p$, Adams operation of the form
\begin{equation}
\label{eq:dual numbers}
\psi^p(a + b\varepsilon) ~=~ a + ba_p \varepsilon
\end{equation}
for some $p$-divisible element $a_p$ in $B$. 

Conversely, let $(a_p)$ be a sequence of elements in $B$ indexed by the set of primes such that $a_p$ is $p$-divisible in $B$.  Then there is a unique isomorphism class of filtered $\lambda$-ring structure over $B \lbrack \varepsilon \rbrack$ with Adams operations satisfying \eqref{eq:dual numbers}.
\end{cor}

%%%%%%%%%%%%%%%%%%%%%%%%%%%%%%%%%%
%%%%%%%%%%%%%%%%%%%%%%%%%%%%%%%%%%

\subsection{Possibly truncated polynomial rings and power series rings}
\label{subsec:polynomial ring}

One might wonder whether or not Theorem \ref{thm:dual numbers} still holds if the ring of dual numbers is replaced by a general truncated polynomial ring $A \lbrack x \rbrack/x^N$ $(N \geq 2)$, a polynomial ring $A \lbrack x \rbrack$, or even a power series ring $A \llbrack x \rrbrack$.   As we will see shortly, there is indeed a corresponding result in each of these cases, but we have to restrict the ground rings to non-zero, $\bZ$-torsionfree $\bQ$-algebras.  As usual, we filter the rings $A \lbrack x \rbrack/x^N$, $A \lbrack x \rbrack$, and $A \llbrack x \rrbrack$ by the $x$-adic filtration.

\bigskip
\begin{thm}
\label{thm:polynomial ring}
Let $A$ be any non-zero, $\bZ$-torsionfree $\bQ$-algebra and let $R$ be one of the following filtered rings:
\begin{itemize}
\item the truncated polynomial ring $A \lbrack x \rbrack/x^N$ $(N \geq 2)$,
\item the polynomial ring $A \lbrack x \rbrack$, or
\item the power series ring $A \llbrack x \rrbrack$. 
\end{itemize}
Then $R$ admits uncountably many mutually non-isomorphic filtered $\lambda$-ring structures.
\end{thm}

\begin{proof}
The proof of this theorem proceeds in exactly the same way as in the proof of Theorem \ref{thm:dual numbers}, so we will not give all the details.  The only notable difference is that the congruence relation
\[
\psi^p(f(x)) \equiv f(x)^p \quad (\text{mod }pR)
\]
for any $f(x)$ in $R$ is now trivially true, since $p$ is a unit in $A$ and hence in $R$.  This is the reason why we need $A$ to be a $\bQ$-algebra and not merely a binomial domain.
\end{proof}

It should be remarked that in Theorem \ref{thm:polynomial ring} for a polynomial ring, the assumption about $\bQ$-algebra cannot be omitted.  Indeed, as shown by Clauwens \cite{clauwens}, the polynomial ring $\bZ \lbrack x \rbrack$ admits essentially only two $\lambda$-ring structures.

%%%%%%%%%%%%%%%%%%%%%%%%%%%%%
%%%%%%%%%%%%%%%%%%%%%%%%%%%%%

\subsection{A Hasse principle for power series $\lambda$-rings}

In this final section we make an observation about distinguishing between non-isomorphic filtered $\lambda$-ring structures over a filtered power series ring.  

In Theorem \ref{thm:polynomial ring} when $R$ is the filtered power series ring $A \llbrack x \rrbrack$, the various filtered $\lambda$-ring structures over it which we constructed are easily distinguishable by considering the leading coefficients in their Adams operations applied to $x$.  However, there are other situations in which one would like to distinguish between different filtered $\lambda$-ring structures over a filtered power series ring where the corresponding Adams operations have the same leading coefficients.  A case in point is the genus of the classifying space $BSU(2)$.  As the author shows in \cite{yau1}, spaces in its genus (which there are uncountably many of them up to homotopy, see \cite{moller}) all have the same $K$-theory filtered ring, $\bZ \llbrack x \rrbrack$ with $x$ in filtration precisely $4$, up to isomorphism.  Therefore, for each one of these $K$-theory filtered $\lambda$-rings, the Adams operation $\psi^p$, $p$ any prime, satisfies the congruence condition
\[
\psi^p(x) ~\equiv~ p^2 x \quad (\text{mod }x^2).
\]
But as Notbohm shows in \cite{notbohm}, these $\lambda$-rings are mutually non-isomorphic.  This raises the question: 
\begin{quote}
\textit{How can filtered $\lambda$-ring structures over a filtered power series ring be distinguished when their corresponding Adams operations have the same leading coefficients?}
\end{quote}

Here we give an answer to this question when the ground ring is a domain (that is, commutative ring with units without zero divisors).  This applies, for example, to the situation above when the filtered ring is $\bZ \llbrack x \rrbrack$.

\bigskip
\begin{thm}
\label{thm:adams operations}
Let $R$ be a domain and let $S_1$ and $S_2$ be two filtered $\lambda$-ring structures over the filtered power series ring $R \llbrack x \rrbrack$ equipped with the $x$-adic filtration, where $x$ is in some fixed positive filtration.  Denote the Adams operations of $S_1$ and $S_2$ by $\psi^n_1$ and $\psi^n_2$, respectively.  Assume that for each prime $p$ there exists an element $\alpha_p \in R$, which is neither $0$ nor a root of unity, such that
\[
\psi^p_1(x) ~\equiv~ \alpha_p x ~\equiv~ \psi^p_2(x) \quad (\text{mod }x^2).
\]
If $\varphi \colon S_1 \to S_2$ is a filtered ring homomorphism $($resp.\ isomorphism$)$, then it is a filtered $\lambda$-ring map $($resp.\ isomorphism$)$ if and only if there exists a prime $p$ such that $\varphi \psi^p_1 = \psi^p_2 \varphi$.
\end{thm}

This Theorem is an immediate consequence of the following observation about commuting power series, which is a slight generalization of a result due to Lubin \cite[Prop.\ 1.1]{lubin}.

\bigskip
\begin{prop}
\label{prop:lubin}
Let $f(x)$ and $g(x)$ be power series over a field $k$.  Suppose that $f(x) \equiv \alpha x \equiv g(x)$ $($mod $x^2)$ for some element $\alpha \in k$ which is neither $0$ nor a root of unity.  Then for every element $c \in k$ there exists a unique power series $h(x)$ over $k$ satisfying the following conditions:
$($1$)$ $h(0) = 0$,
$($2$)$ $h^\prime(0) = c$, and
$($3$)$ $h(g(x)) = f(h(x))$ as power series.
\end{prop}

Accepting this proposition for the moment, let us prove the theorem.

\begin{proof}[Proof of Theorem \ref{thm:adams operations}]
Since the ``only if'' part is obvious, we now consider the ``if'' part.

Note that for any prime $q$ the power series $\psi^q_1(\varphi(x))$ and $\varphi(\psi^q_2(x))$ have the same linear coefficient, both have $0$ constant term, and both satisfy the last condition in Proposition \ref{prop:lubin} with $g = \psi^p_1$ and $f = \psi^p_2$.  Therefore, by uniqueness we conclude that $\psi^q_1 (\varphi(x))$ and $\varphi (\psi^q_2(x))$ coincide as power series, and so $\varphi \psi^q_1$ is equal to $\psi^q_2 \varphi$ as filtered ring endomorphisms on $R \llbrack x \rrbrack$.

This proves the theorem.
\end{proof}

Since the proof of Proposition \ref{prop:lubin} is an elementary degree-by-degree argument, which is almost exactly the same as that of Lubin's \cite[Prop.\ 1.1]{lubin}, we will not give the details.  The only important point in the argument is that one can divide by $(\alpha^j - \alpha)$ in $k$ for any integer $j \geq 2$, which is possible by the hypothesis on $\alpha$.

                                   %%  References  %%

\end{document}